\documentclass[11pt]{amsart}

\usepackage{amsmath,amsthm,amscd,euscript,longtable}

\newtheorem{Th}{\hskip\parindent Theorem}

\newcommand{\R}{\mathcal{R}}
\newcommand{\K}{\mathcal{K}}

\usepackage{tikz}

\title{Supercontinuants}
\author{Alexey Ustinov}
\date{05 March 2015}

\thanks{The research was performed under support of Russian Science Foundation (Project N~14-11-00335)}


\email[Alexey Ustinov]{ustinov@iam.khv.ru}

\address{Institute of Applied Mathematics,
                                        Khabarovsk Division, Russian Academy of Sciences,\\ 54 Dzerzhinsky Street, Khabarovsk, 680000, Russia}

\begin{document}

\begin{abstract} Morier-Genoud, Ovsienko and  Tabachnikov introduced supersymmetric frieze patterns (see ArXiv 1501.07476). This note gives a solution to Problem 1 from that article: determine the formula for the entries of a superfrieze.\end{abstract}

\maketitle

This note gives a solution to  Problem 1 from the article~\cite{Morier-Genoud2015}: determine the formula for the entries of a superfrieze.

Let $\R=\R_0\oplus\R_1$ be an arbitrary supercommutative
ring, and the sequences $\{v_i\}$, $\{w_i\}$, with $v_i\in\R_0$, $w_i\in\R_1$, be defined by the initial conditions $v_{-1}=0$, $v_0=1$, $w_0=0$ and the recurrence relation
\begin{equation}
\label{1}
v_i=a_iv_{i-1}-v_{i-2}-\beta_iw_{i-1},\quad w_i=w_{i-1}+\beta_iv_{i-1}\quad(i\in \mathbb{Z}).
\end{equation}
In particular,
\begin{gather*}
v_1=a_1,\quad v_2=a_1a_2-1+\beta_1\beta_2,\quad v_3=a_1a_2a_3-a_1-a_3+a_1\beta_2\beta_3+a_3\beta_1\beta_2+\beta_1\beta_3;\\
w_1=\beta_1,\quad w_2=a_1\beta_2+\beta_1,\quad w_3= a_1a_2\beta_3+a_1\beta_2+\beta_1\beta_2\beta_3+\beta_1-\beta_3.
\end{gather*}
The problem is to express $v_n$, $w_n$ in terms of $a_1$, \ldots, $a_n$ and $\beta_1$, \ldots, $\beta_n$. Such expression will be called  \textit{supercontinuants}. (See~\cite{Knu} for the properties of the classical continuants.)

We define two sequences of supercontinuants $$\{\K\bigl(\begin{smallmatrix} a_1\\
\begin{smallmatrix} \beta_1 & \beta_1
\end{smallmatrix}
\end{smallmatrix}|\begin{smallmatrix}
a_2\\\begin{smallmatrix} \beta_2 & \beta_2
\end{smallmatrix}
\end{smallmatrix}|\ldots |\begin{smallmatrix}
a_n\\\begin{smallmatrix} \beta_n & \beta_n
\end{smallmatrix}
\end{smallmatrix}\bigr)\}\text{\quad and \quad} \{\K\bigl(\begin{smallmatrix}  a_1\\\begin{smallmatrix}
\beta_1 & \beta_1
\end{smallmatrix}
\end{smallmatrix}|\ldots |\begin{smallmatrix}
a_{n-1}\\\begin{smallmatrix} \beta_{n-1} & \beta_{n-1}
\end{smallmatrix}
\end{smallmatrix}|\beta_n\bigr)\}$$ by the initial conditions  $\K()=1$, $\K(\begin{smallmatrix} a_1\\
\begin{smallmatrix} \beta_1 & \beta_1
\end{smallmatrix}
\end{smallmatrix})=a_1$, $\K(\beta_1)=\beta_1$  and the recurrence relations
\begin{gather}\nonumber
\K\bigl(\begin{smallmatrix} a_1\\
	\begin{smallmatrix} \beta_1 & \beta_1
	\end{smallmatrix}
\end{smallmatrix}|\ldots |\begin{smallmatrix}
a_n\\\begin{smallmatrix} \beta_n & \beta_n
\end{smallmatrix}
\end{smallmatrix}\bigr)=a_n\K\bigl(\begin{smallmatrix} a_1\\
\begin{smallmatrix} \beta_1 & \beta_1
\end{smallmatrix}
\end{smallmatrix}|\ldots |\begin{smallmatrix}
a_{n-1}\\\begin{smallmatrix} \beta_{n-1} & \beta_{n-1}
\end{smallmatrix}
\end{smallmatrix}\bigr)-\K\bigl(\begin{smallmatrix} a_1\\
\begin{smallmatrix} \beta_1 & \beta_1
\end{smallmatrix}
\end{smallmatrix}|\ldots |\begin{smallmatrix}
a_{n-2}\\\begin{smallmatrix} \beta_{n-2} & \beta_{n-2}
\end{smallmatrix}
\end{smallmatrix}\bigr)
\\
\label{2}-\beta_n\K\bigl(\begin{smallmatrix} a_1\\
	\begin{smallmatrix} \beta_1 & \beta_1
	\end{smallmatrix}
\end{smallmatrix}|\ldots |\begin{smallmatrix}
a_{n-2}\\\begin{smallmatrix} \beta_{n-2} & \beta_{n-2}
\end{smallmatrix}
\end{smallmatrix}|\beta_{n-1}\bigr),\\
\nonumber \K\bigl(\begin{smallmatrix}  a_1\\\begin{smallmatrix}
		\beta_1 & \beta_1
	\end{smallmatrix}
\end{smallmatrix}|\ldots |\begin{smallmatrix}
a_{n-1}\\\begin{smallmatrix} \beta_{n-1} & \beta_{n-1}
\end{smallmatrix}
\end{smallmatrix}|\beta_n\bigr)=\beta_n\K\bigl(\begin{smallmatrix}  a_1\\\begin{smallmatrix}
	\beta_1 & \beta_1
\end{smallmatrix}
\end{smallmatrix}|\ldots |\begin{smallmatrix}
a_{n-1}\\\begin{smallmatrix} \beta_{n-1} & \beta_{n-1}
\end{smallmatrix}
\end{smallmatrix}\bigr)+\K\bigl(\begin{smallmatrix}  a_1\\\begin{smallmatrix}
	\beta_1 & \beta_1
\end{smallmatrix}
\end{smallmatrix}|\ldots |\begin{smallmatrix}
a_{n-2}\\\begin{smallmatrix} \beta_{n-2} & \beta_{n-2}
\end{smallmatrix}
\end{smallmatrix}|\beta_{n-1}\bigr).
\end{gather}
%(We assume that, for each couple of arguments
%$(\beta_i,\beta_i)$, the element $a_i$ is also included in
%the list of arguments of
%$\K(\ldots;\beta_i,\beta_i;\ldots)$, and we suppress $a_i$
%from the notation.) 
From~\eqref{1} and~\eqref{2} it easily
follows that
\begin{equation*}
v_n=\K\bigl(\begin{smallmatrix} a_1\\
	\begin{smallmatrix} \beta_1 & \beta_1
	\end{smallmatrix}
\end{smallmatrix}|\ldots |\begin{smallmatrix}
a_n\\\begin{smallmatrix} \beta_n & \beta_n
\end{smallmatrix}
\end{smallmatrix}\bigr), \qquad
w_n=\K\bigl(\begin{smallmatrix}  a_1\\\begin{smallmatrix}
		\beta_1 & \beta_1
	\end{smallmatrix}
\end{smallmatrix}|\ldots |\begin{smallmatrix}
a_{n-1}\\\begin{smallmatrix} \beta_{n-1} & \beta_{n-1}
\end{smallmatrix}
\end{smallmatrix}|\beta_n\bigr).
\end{equation*}

The classical continuants $K(a_1,\ldots,a_n)$,
corresponding to reduced regular continued fractions
$$a_1-\cfrac{1}{a_2-{\atop\ddots\,\displaystyle{-\cfrac{1}{a_n}}}},$$ are defined by
$$K()=1,\quad K(a_1)=a_1,\quad K(a_1,\ldots,a_n)=a_nK(a_1,\ldots,a_{n-1})-K(a_1,\ldots,a_{n-2}).$$
There is Euler's rule which allows one to write down all summands of $K(a_1,\ldots,a_n)$: starting with
the product $a_1 a_2 \ldots a_n$, we strike out adjacent pairs $a_ia_{i+1}$ in all
possible ways. If a pair $a_ia_{i+1}$ is struck out, then it must be replaced by $-1$. We can represent Euler's rule graphically by constructing all
``Morse code'' sequences of dots and dashes having length $n$, where each dot
contributes $1$ to the length and each dash contributes $2$. For example $K(a_1,a_2,a_3,a_4)$ consists of the following summands:
\begin{center}
                                        \begin{tikzpicture}[scale=.6]
                                        %\draw  (0,0) -- (1,0);
                                        \draw  (2,0) -- (3,0);\draw  (0,1) -- (1,1);

                                        \filldraw (0,0) circle (0.07);\filldraw (0,1) circle (0.07);\filldraw (0,2) circle (0.07);
                                        \filldraw (1,0) circle (0.07);\filldraw (1,1) circle (0.07);\filldraw (1,2) circle (0.07);
                                        \filldraw (2,0) circle (0.07);\filldraw (2,1) circle (0.07);\filldraw (2,2) circle (0.07);
                                        \filldraw (3,0) circle (0.07);\filldraw (3,1) circle (0.07);\filldraw (3,2) circle (0.07);

                                        \filldraw [right] (5,1) node {$\mapsto -a_3a_4$};
                                        \filldraw [right] (5,0) node {$\mapsto -a_1a_2$};
                                        \filldraw [right] (5,2) node {$\mapsto a_1a_2a_3a_4$};

                                        \begin{scope}[scale=1, xshift=280]

                                        \filldraw (0,1) circle (0.07);\filldraw (0,2) circle (0.07);
                                        \filldraw (1,1) circle (0.07);\filldraw (1,2) circle (0.07);
                                        \filldraw (2,1) circle (0.07);\filldraw (2,2) circle (0.07);
                                        \filldraw (3,1) circle (0.07);\filldraw (3,2) circle (0.07);

                                        \draw  (0,1) -- (1,1);\draw  (2,1) -- (3,1);\draw  (2,2) -- (1,2);

                                        \filldraw [right] (5,1) node {$\mapsto 1$};
                                        \filldraw [right] (5,2) node {$\mapsto -a_1a_4$};

                                        \end{scope}
                                        \end{tikzpicture}
\end{center}

By analogy with  Euler's rule, we can construct a similar rule for calculation of supercontinuants.

\begin{Th}
The summands of $\K\bigl(\begin{smallmatrix} a_1\\
\begin{smallmatrix} \beta_1 & \beta_1
\end{smallmatrix}
\end{smallmatrix}|\begin{smallmatrix}
a_2\\\begin{smallmatrix} \beta_2 & \beta_2
\end{smallmatrix}
\end{smallmatrix}|\ldots \bigr)$ can be obtained from the product $\beta_1\beta_1\beta_2\beta_2\ldots$ by the following rule:  we strike out adjacent pairs and adjacent 4-tuples $\beta_i\beta_i\beta_{i+1}\beta_{i+1}$ in all
possible ways; for deleted  pairs  and 4-tuples we make the substitutions
$\beta_i\beta_i\to a_i$, $\beta_i \beta_{i+1}\to 1$, $\beta_i\beta_i\beta_{i+1}\beta_{i+1}\to -1.$
\end{Th}

This rule can be represented graphically as well. To each monomial there corresponds a sequence of  total length $2n$ (or $2n-1$) consisting of dots (of the length one), dashes (of the length two) and long dashes (of the length four). For example, the monomials of $\K\bigl(\begin{smallmatrix} a_1\\
\begin{smallmatrix} \beta_1 & \beta_1
\end{smallmatrix}
\end{smallmatrix}|\begin{smallmatrix}
a_2\\\begin{smallmatrix} \beta_2 & \beta_2
\end{smallmatrix}
\end{smallmatrix}|\beta_{3}\bigr)$ can be obtained from the product $\beta_1\beta_1 \beta_2\beta_2\beta_{3}$ as follows:
\begin{center}
\begin{tikzpicture}[scale=.6]
\draw  (0,1) -- (1,1);\draw  (3,1) -- (4,1);\draw  (0,2) -- (1,2);\draw  (2,2) -- (3,2);

\filldraw (0,0) circle (0.07);\filldraw (0,1) circle (0.07);\filldraw (0,2) circle (0.07);
\filldraw (1,0) circle (0.07);\filldraw (1,1) circle (0.07);\filldraw (1,2) circle (0.07);
\filldraw (2,0) circle (0.07);\filldraw (2,1) circle (0.07);\filldraw (2,2) circle (0.07);
\filldraw (3,0) circle (0.07);\filldraw (3,1) circle (0.07);\filldraw (3,2) circle (0.07);
\filldraw (4,0) circle (0.07);\filldraw (4,1) circle (0.07);\filldraw (4,2) circle (0.07);

\filldraw [right] (5,1) node {$\mapsto a_1\beta_2$}; \draw  (1,0) -- (2,0);    \filldraw [right] (5,0) node {$\mapsto\beta_1\beta_2\beta_3$};
\filldraw [right] (5,2) node {$\mapsto a_1a_2\beta_3$};

\begin{scope}[scale=1, xshift=280]

\filldraw (0,1) circle (0.07);\filldraw (0,2) circle (0.07);
\filldraw (1,1) circle (0.07);\filldraw (1,2) circle (0.07);
\filldraw (2,1) circle (0.07);\filldraw (2,2) circle (0.07);
\filldraw (3,1) circle (0.07);\filldraw (3,2) circle (0.07);
\filldraw (4,1) circle (0.07);\filldraw (4,2) circle (0.07);

\draw  (1,2) -- (2,2);\draw  (3,2) -- (4,2);\draw  (0,1) -- (3,1);

\filldraw [right] (5,1) node {$\mapsto-\beta_3$};
\filldraw [right] (5,2) node {$\mapsto\beta_1$};

\end{scope}
\end{tikzpicture}
\end{center}
Let us note that the odd
variables anticommute with each other. In particular, $\beta_i^2=0$, and in each pair $\beta_i\beta_i$ at least one variable must be struck out. Supercontinuants become the usual continuants if all odd variables are replaced by zeros.

Supercontinuants can be expressed as determinants.

\begin{Th}
\begin{gather}\label{3}
                                   \K\bigl(\begin{smallmatrix} a_1\\
                                   	\begin{smallmatrix} \beta_1 & \beta_1
                                   	\end{smallmatrix}
                                   \end{smallmatrix}|\ldots |\begin{smallmatrix}
                                   a_n\\\begin{smallmatrix} \beta_n & \beta_n
                                   \end{smallmatrix}
                                \end{smallmatrix}\bigr)=\begin{vmatrix}
                                        %  \begin{array}{cccccc}
                                        a_1 & -1+\beta_1\beta_2 & \beta_1\beta_3 & \cdots & \beta_1\beta_{n-1} & \beta_1\beta_n \\
                                        -1 & a_2 & -1+\beta_2\beta_3 & \cdots & \beta_2\beta_{n-1} & \beta_2\beta_n \\
                                        0 & -1 & a_3 & \cdots  & \beta_3\beta_{n-1} & \beta_3\beta_n  \\
                                        \cdots & \cdots  & \cdots & \cdots & \cdots & \cdots \\
                                        0 & \cdots & 0 & -1 & a_{n-1} & -1+\beta_{n-1}\beta_n \\
                                        0 & 0  & \cdots & 0 & -1 & a_n \\
                                        % \end{array}
                                        \end{vmatrix},\\\nonumber
                                        \K\bigl(\begin{smallmatrix}  a_1\\\begin{smallmatrix}
                                        		\beta_1 & \beta_1
                                        	\end{smallmatrix}
                                        \end{smallmatrix}|\ldots |\begin{smallmatrix}
                                        a_{n-1}\\\begin{smallmatrix} \beta_{n-1} & \beta_{n-1}
                                        \end{smallmatrix}
                                    \end{smallmatrix}|\beta_n\bigr)=\begin{vmatrix}
                                        %  \begin{array}{cccccc}
                                        a_1 & -1+\beta_1\beta_2 & \beta_1\beta_3 & \cdots & \beta_1\beta_{n-1} & \beta_1\\
                                        -1 & a_2 & -1+\beta_2\beta_3 & \cdots & \beta_2\beta_{n-1} & \beta_2 \\
                                        0 & -1 & a_3 & \cdots  & \beta_3\beta_{n-1} & \beta_3\\
                                        \cdots & \cdots  & \cdots & \cdots & \cdots & \cdots \\
                                        0 & \cdots & -1 & a_{n-2} & -1+\beta_{n-2}\beta_{n-1}  & \beta_{n-2}\\
                                        0 & \cdots & 0 & -1 & a_{n-1} & \beta_{n-1}\\
                                        0 & 0  & \cdots & 0 & -1 & \beta_n \\
                                        % \end{array}
                                        \end{vmatrix}.
                                        \end{gather}
\end{Th}

The second determinant in  Theorem 2 is well-defined
because odd variables occupy only one column. The proofs of
Theorems 1 and 2 follow by induction from recurrence
relations~\eqref{2}, and we do not dwell on them.

The supercontinuants of the form $
\K\bigl(\beta_1|\begin{smallmatrix}  a_2\\\begin{smallmatrix}
\beta_2 & \beta_2
\end{smallmatrix}
\end{smallmatrix}|\ldots |\begin{smallmatrix}
a_{n-1}\\\begin{smallmatrix} \beta_{n-1} & \beta_{n-1}
\end{smallmatrix}
\end{smallmatrix}|\beta_n\bigr)
$ also may be defined by the rule from the Theorem 1. For example
$$\K(\beta_1|\beta_2)=\beta_1\beta_2+1,\quad \K\bigl(\beta_1|\begin{smallmatrix}  a_2\\\begin{smallmatrix}
\beta_2 & \beta_2
\end{smallmatrix}
\end{smallmatrix}|\beta_3\bigr)=a_2\beta_1\beta_3+\beta_1\beta_2+\beta_2\beta_3+1.$$ These supercontinuants can be represented in terms of determinants as well (we assume that the determinant is expanded in the first column, and the same rule is applied to all determinants of  smaller matrices).

\begin{Th}
The supercontinuants $\K\bigl(\beta_1|\begin{smallmatrix}  a_2\\\begin{smallmatrix}
\beta_2 & \beta_2
\end{smallmatrix}
\end{smallmatrix}|\ldots |\begin{smallmatrix}
a_{n-1}\\\begin{smallmatrix} \beta_{n-1} & \beta_{n-1}
\end{smallmatrix}
\end{smallmatrix}|\beta_n\bigr)$ satisfy the
recurrence relation
\begin{equation}
\label{5} \K\bigl(\beta_1|\begin{smallmatrix}  a_2\\\begin{smallmatrix}
		\beta_2 & \beta_2
	\end{smallmatrix}
\end{smallmatrix}|\ldots |\begin{smallmatrix}
a_{n-1}\\\begin{smallmatrix} \beta_{n-1} & \beta_{n-1}
\end{smallmatrix}
\end{smallmatrix}|\beta_n\bigr)=-\beta_n\K\bigl(\beta_1|\begin{smallmatrix}  a_2\\\begin{smallmatrix}
\beta_2 & \beta_2
\end{smallmatrix}
\end{smallmatrix}|\ldots |\begin{smallmatrix}
a_{n-1}\\\begin{smallmatrix} \beta_{n-1} & \beta_{n-1}
\end{smallmatrix}
\end{smallmatrix}\bigr)+\K\bigl(\beta_1|\begin{smallmatrix}  a_2\\\begin{smallmatrix}
	\beta_2 & \beta_2
\end{smallmatrix}
\end{smallmatrix}|\ldots |\beta_{n-1}\bigr)
\qquad(n\ge 2)
\end{equation}
and can be expressed in the following form:
\begin{gather*}
\K\bigl(\beta_1|\begin{smallmatrix}  a_2\\\begin{smallmatrix}
		\beta_2 & \beta_2
	\end{smallmatrix}
\end{smallmatrix}|\ldots |\begin{smallmatrix}
a_{n-1}\\\begin{smallmatrix} \beta_{n-1} & \beta_{n-1}
\end{smallmatrix}
\end{smallmatrix}|\beta_n\bigr)=\begin{vmatrix}
                                        %  \begin{array}{cccccc}
                                        \beta_1 & \beta_2 & \beta_3 & \cdots & \beta_{n-1} & 1\\
                                        -1 & a_2 & -1+\beta_2\beta_3 & \cdots & \beta_2\beta_{n-1} & \beta_2 \\
                                        0 & -1 & a_3 & \cdots  & \beta_3\beta_{n-1} & \beta_3\\
                                        \cdots &  \cdots & \cdots & \cdots & \cdots & \cdots \\
                                        0 & \cdots & -1 & a_{n-2} & -1+\beta_{n-2}\beta_{n-1}  & \beta_{n-2}\\
                                        0 & \cdots & 0 & -1 & a_{n-1} & \beta_{n-1}\\
                                        0 & 0  & \cdots & 0 & -1 & \beta_n \\
                                        % \end{array}
                                        \end{vmatrix}.
                                        \end{gather*}
\end{Th}

The proof of  formula~\eqref{5} is an application of the
rule from Theorem 1. The determinant formula follows by
induction from the recurrence relation~\eqref{5}.

Finally, the even supercontinuants  $\K\bigl(\begin{smallmatrix} a_1\\
\begin{smallmatrix} \beta_1 & \beta_1
\end{smallmatrix}
\end{smallmatrix}|\ldots |\begin{smallmatrix}
a_n\\\begin{smallmatrix} \beta_n & \beta_n
\end{smallmatrix}
\end{smallmatrix}\bigr)$  can be also expressed as Berezinians. Recall that the Berezinian of the matrix
$$
\mathrm{Ber}
\begin{pmatrix}
A&B\\C&D\\
\end{pmatrix},
$$
where where $A$ and $D$ have even entries, and $B$ and $C$
have odd entries, is given by the formula
\begin{equation}
\label{Ber} \det(A-BD^{-1}C)\det(D)^{-1},
\end{equation}see, e.g., \cite{Be}.

\begin{Th}
\begin{gather*}
 \K\bigl(\begin{smallmatrix} a_1\\
 	\begin{smallmatrix} \beta_1 & \beta_1
 	\end{smallmatrix}
 \end{smallmatrix}|\ldots |\begin{smallmatrix}
 a_n\\\begin{smallmatrix} \beta_n & \beta_n
 \end{smallmatrix}
\end{smallmatrix}\bigr)=\mathrm{Ber}\begin{pmatrix}
A&B\\C&D\\                              \end{pmatrix},
\end{gather*}
where
\begin{gather*}
A=\begin{pmatrix}
a_1&-1&0&\cdots&0\\-1&a_2&-1&\ddots&\vdots\\
0&-1&a_3&\ddots&0\\
\vdots&\ddots&\ddots&\ddots&-1\\
0&\cdots&0&-1&a_n\\
\end{pmatrix},\qquad B=\begin{pmatrix}
\beta_1&\beta_2 &\beta_3&\cdots&\beta_n \\0&\beta_2&\beta_3&\ddots&\vdots\\
0&0&\beta_3&\ddots&\beta_n\\
\vdots&\ddots&\ddots&\ddots&\beta_n\\
0&\cdots&0&0&\beta_n\\
\end{pmatrix},\\
C=\begin{pmatrix}
-\beta_1&0&0&\cdots&0\\0&-\beta_2&0&\ddots&\vdots\\
0&0&-\beta_3&\ddots&0\\
\vdots&\ddots&\ddots&\ddots&0\\
0&\cdots&0&0&-\beta_n\\
\end{pmatrix},\qquad D=\begin{pmatrix}
1&0 &0&\cdots&0 \\0&1&0&\ddots&\vdots\\
0&0&1&\ddots&0\\
\vdots&\ddots&\ddots&\ddots&0\\
0&\cdots&0&0&1\\
\end{pmatrix}.
\end{gather*}
\end{Th}

Theorem 4 is direct corollary of~\eqref{3}
and~\eqref{Ber}.

It follows from recurrence relations~\eqref{2} and~\eqref{5} that  the
number of terms in supercontinuants $ \K\bigl(\begin{smallmatrix} a_1\\
\begin{smallmatrix} \beta_1 & \beta_1
\end{smallmatrix}
\end{smallmatrix}|\ldots |\begin{smallmatrix}
a_n\\\begin{smallmatrix} \beta_n & \beta_n
\end{smallmatrix}
\end{smallmatrix}\bigr)$, $\K\bigl(\begin{smallmatrix}  a_1\\\begin{smallmatrix}
\beta_1 & \beta_1
\end{smallmatrix}
\end{smallmatrix}|\ldots |\begin{smallmatrix}
a_{n-1}\\\begin{smallmatrix} \beta_{n-1} & \beta_{n-1}
\end{smallmatrix}
\end{smallmatrix}|\beta_n\bigr)$ and $\K\bigl(\beta_1|\begin{smallmatrix}  a_2\\\begin{smallmatrix}
\beta_2 & \beta_2
\end{smallmatrix}
\end{smallmatrix}|\ldots |\begin{smallmatrix}
a_{n-1}\\\begin{smallmatrix} \beta_{n-1} & \beta_{n-1}
\end{smallmatrix}
\end{smallmatrix}|\beta_n\bigr)$
coincide respectively with the sequences (see~\cite{OEIS})
\begin{align*}
A077998&:1, 3, 6, 14, 31, 70, 157, 353, 793, 1782, 4004,  \ldots\\
A006054&:1, 2, 5, 11, 25, 56, 126, 283, 636, 1429, 3211, \ldots\\
A052534&:1, 2, 4,  9, 20, 45, 101, 227, 510, 1146, 2575,
\ldots
\end{align*}

%{\bf A.\,V.~Ustinov}

%{\it E-mail}: ustinov@iam.khv.ru

\end{document}